\documentclass{article}
\usepackage{amsfonts}
\usepackage{graphicx}
\usepackage{epstopdf}
\usepackage{color,url,bm}
\usepackage{marginnote}
\usepackage{hyperref}
\usepackage{bm,color}

\def\blue{\color[rgb]{0,0,1}}

\newtheorem{rem}{Remark}

\def\RR{\mathbb{R}}

\def\dd{\mathrm{d}}
\def\pmatrix{\left(\begin{array}}
\def\endpmatrix{\end{array}\right)}
\def\bfb{\bm{b}}
\def\bfc{\bm{c}}
\def\bfe{\bm{e}}

\def\bfgamma{\bm{\gamma}}

\def\P{{\cal P}}
\def\I{{\cal I}}

\def\no{\noindent}

\begin{document}
\title{Continuous-Stage Runge-Kutta approximation to Differential Problems}

% Author Orchid ID: enter ID or remove command

% Authors, for the paper (add full first names)
\author{Pierluigi Amodio\thanks{Dipartimento di Matematica, Universit\`a di Bari, Via Orabona 4, 70125 Bari, Italy.} \and Luigi Brugnano\thanks{Dipartimento di Matematica e Informatica ``U.\,Dini'', Universit\`a di Firenze, Viale Morgagni 67/A, 50134~Firenze, Italy.} \and Felice Iavernaro$^*$}

\maketitle

\begin{abstract}In recent years, the efficient numerical solution of Hamiltonian problems has led to the definition of a class of energy-conserving Runge-Kutta methods named {\em Hamiltonian Boundary Value Methods (HBVMs)}. Such methods admit an interesting interpretation in terms of {\em continuous-stage Runge-Kutta methods}, which is here recalled and revisited for general differential problems.  

\bigskip
{\no\bf Keywords:} Hamiltonian problems; Hamiltonian Boundary Value Methods; HBVMs; ODE-IVPs; continuous-stage Runge-Kutta methods; continuous-stage Runge-Kutta-Nystr\"om methods.

\bigskip
{\no\bf MSC:} 65L06; 65P10; 65L05.
\end{abstract}

%\begin{document}

\section{Introduction}
{\em Continuous-stage Runge-Kutta methods (csRK, hereafter)}, introduced by But\-cher \cite{Bu1,Bu2,Bu3}, have been used in recent years as a useful tool for studying energy-conserving methods for Hamiltonian problems  (see, e.g., \cite{BIT2010,BIT2012_3,H2010,MB2016,QMcL2008,T2018,TC2007,TS2012,TS2014}). In particular,
we shall at first consider methods, within this latter class, having Butcher tableau in the form:
\begin{equation}\label{csRK}
\begin{array}{c|c} c & a(c,\tau)\\ \hline & b(c)\end{array}, 
\end{equation}
where $c,\tau\in[0,1]$, and
\begin{equation}\label{actau}
a(c,\tau):[0,1]\times [0,1]\rightarrow\RR, \qquad b(c):[0,1]\rightarrow 1,
\end{equation}
are suitable functions defining the method. Hereafter, we shall use the notation (\ref{actau}), in place of the more commonly used $a_{c\tau},\,b_c$, in order to make clear that these are functions of the respective arguments. For later use, we shall also denote
\begin{equation}\label{a1c}
\dot a(c,\tau) = \frac{\dd}{\dd c} a(c,\tau).
\end{equation}
When used for solving the ODE-IVP\,\footnote{In the sequel, we shall assume $f$ to be analytical.}
\begin{equation}\label{odeivp}
\dot y(t) = f(y(t)), \qquad t\in[0,h], \qquad y(0)=y_0\in\RR^m,
\end{equation}
the method (\ref{csRK})-(\ref{a1c}) provides an approximation
\begin{equation}\label{apprivp}
\dot\sigma(ch) = \int_0^1 \dot a(c,\tau)f(\sigma(\tau h))\dd\tau, \qquad c\in[0,1], \qquad \sigma(0)=y_0,
\end{equation}
to (\ref{odeivp}) and, therefore,
\begin{equation}\label{u}
\sigma(ch) = y_0+h\int_0^1 a(c,\tau)f(\sigma(\tau h))\dd\tau, \qquad c\in[0,1],
\end{equation}
with the approximation to $y(h)$ given by
\begin{equation}\label{y1}
y_1 = y_0 + h\int_0^1 b(c) f(\sigma(ch))\dd c.
\end{equation}
In particular, we shall hereafter consider the natural choiche
\begin{equation}\label{bc}
b(c)\equiv 1, \qquad c\in[0,1],
\end{equation}
though different choices have been also considered \cite{MB2016}.\footnote{Clearly, because of consistenty, one requires that $\int_0^1b(c)\dd c=1$.} The arguments studied in this paper strictly follows those in \cite{ABI2019} (in turn, inspired by \cite{T2018,TSZ2019,TZ2018}), derived from the energy-conserving methods called {\em Hamiltonian Boundary Value Methods (HBVMs)}, which have been the subject of many investigations \cite{ABI2015, ABI2019, ABI2020, ABI2022, BBFCI2018, BBTZ2020,BCMR2012, BFCI2013, BFC2014, BFCI2014, BFCI2015, BFCI2019, BGIW2018, BGZ2019, LIMbook2016, BI2018ax, BIMR2019, BIS2009, BIS2010, BI2012, BIT2009, BIT2010, BIT2011, BIT2012, BIT2012_2, BIT2012_3, BIT2015, BIZ2020,  BMR2019, BMR2019_1, BGS2019, BrSu2014, BZL2018}. Here, we consider also a relevant generalization, w.r.t. the arguments in \cite{ABI2019}, as explained below. With this premise, the structure of the paper is as follows: in Section~\ref{ode1} we recall the basic facts concerning the initial value problem (\ref{odeivp}); in Section~\ref{ode2} we generalize the approach to the case of second-order problems, thus obtaining a {\em continuous-stage Runge-Kutta-Nystr\"om (csRKN, hereafter) method}, whose generalization for $k$th-order problems is also sketched; Section~\ref{discre} is devoted to derive relevant families of methods for the previous classes of problems; at last, in Section~\ref{fine} a few conclusions are given.

\section{Approximation of ODE-IVPs}\label{ode1}
Let us consider, at first, the initial value problem (\ref{odeivp}). As done in \cite{ABI2019}, for our analysis we shall use an expansion of the vector field along the orthonormal Legendre polynomial basis:
\begin{equation}\label{orto}
P_i\in\Pi_i, \qquad \int_0^1 P_i(x)P_j(x)\dd x =\delta_{ij}, \qquad i,j=0,1,\dots,
\end{equation}
where, as usual, $\Pi_i$ denotes the vector space of polynomials of degree $i$, and $\delta_{ij}$ is the Kronecker symbol.
Consequently, we can rewrite (\ref{odeivp}) as
\begin{equation}\label{exp1}
\dot y(ch) = \sum_{j\ge0} P_j(c)\gamma_j(y), \qquad c\in[0,1],\qquad y(0) = y_0\in\RR^m,
\end{equation}
with
\begin{equation}\label{gammaj}
\gamma_j(y) = \int_0^1 P_j(\tau) f(y(\tau h))\dd\tau, \qquad j=0,1,\dots.
\end{equation}
Integrating side-by-side, and imposing the initial condition, one then obtains that the solution of (\ref{exp1}) is formally given by
\begin{equation}\label{sol1}
y(ch) = y_0+h\sum_{j\ge0} \int_0^c P_j(x)\dd x\,\gamma_j(y), \qquad c\in[0,1].
\end{equation}
For $c=1$, by considering that, by virtue of (\ref{orto}), $\int_0^1P_j(x)\dd x=\delta_{j0}$, and taking into account (\ref{odeivp}), one obtains:
\begin{equation}\label{yh}
y(h) = y_0 + h\int_0^1 f(y(ch))\dd c \,\equiv\, y_0+\int_0^h \dot y(t)\dd t,
\end{equation}
i.e., the Fundamental Theorem of the Calculus. Interestingly, by setting
\begin{equation}\label{ainf}
a_\infty(c,\tau) = \sum_{j\ge0} \int_0^c P_j(x)\dd x P_j(\tau),
\end{equation}
one obtains that (\ref{exp1})-(\ref{gammaj}) and (\ref{sol1})-(\ref{yh}) can be rewritten, respectively, as:
\begin{eqnarray}\nonumber
\dot y(ch) &=& \int_0^1\dot a_\infty(c,\tau)f(y(\tau h))\dd\tau, \qquad y(0)=y_0\in\RR^m,\\ \label{csRKinf}
y(ch)        &=& y_0 + h\int_0^1 a_\infty(c,\tau)f(y(\tau h))\dd\tau, \qquad c\in[0,1],\\ \nonumber
y(h)          &=& y_0 + h\int_0^1 f(y(ch))\dd c.
\end{eqnarray} 
Consequently, 
\begin{equation}\label{csRK1}
\begin{array}{c|c} c & a_\infty(c,\tau)\\ \hline & 1\end{array}, 
\end{equation}
is the csRK ``method'' providing the exact solution to the problem.

\subsection{Vector formulation}\label{ode1_v}

For later use, we now cast the formulation of (\ref{ainf}) in vector form. For this purpose, let us introduce the infinite vectors,
\begin{equation}\label{PIinf}
\P_\infty(c) = \pmatrix{c} P_0(c)\\ P_1(c)\\ \vdots\endpmatrix, \qquad \I_\infty(c) = \int_0^c \P_\infty(x)\dd x,
\end{equation}
 and the matrix,
 \begin{equation}\label{Xinf}
 X_\infty = \pmatrix{cccc} \xi_0 &-\xi_1\\
                                         \xi_1 & 0 & -\xi_2\\
                                                  &\xi_2 &\ddots &\ddots \\
                                                  &        &\ddots
                                                  \endpmatrix,
                                                  \qquad \xi_i = \frac{1}{2\sqrt{|4i^2-1|}}, \qquad i=0,1,\dots,
                                              \end{equation}    
also recalling that, by virtue of (\ref{orto}), and due to the well-known relations between the Legendre polynomials and their integrals,
\begin{equation}\label{PPPIXinf} 
\I_\infty(c)^\top = P_\infty(c)^\top X_\infty, ~
\int_0^1 \P_\infty(\tau)\P_\infty(\tau)^\top\dd \tau = I, ~ \int_0^1 \P_\infty(\tau)\I_\infty(\tau)^\top\dd \tau = X_\infty,
\end{equation}
with $I$ the identity operator. Consequently, with reference to (\ref{ainf}), one has
\begin{equation}\label{ainfeq}
a_\infty(c,\tau) = \I_\infty(c)^\top\P_\infty(\tau) = \P_\infty(c)^\top X_\infty \P_\infty(\tau),
\end{equation}                                              
and, in particular, one may regard the Butcher tableau, equivalent to (\ref{csRK1})
\begin{equation}\label{csRK1W}
\begin{array}{c|c} c & \P_\infty(c)^\top X_\infty \P_\infty(\tau) \\ \hline & 1\end{array}, 
\end{equation}
as the corresponding $W$-transformation \cite{HW2002} of the continuous problem.

It is worth mentioning that, by defining the infinite vector (see (\ref{gammaj})), 
\begin{equation}\label{bgam}
\bfgamma := \pmatrix{c} \gamma_0(y)\\ \gamma_1(y)\\ \vdots\endpmatrix \equiv \int_0^1 \P_\infty(\tau)\otimes I_m f(y(\tau h))\dd \tau,
\end{equation}
then (\ref{sol1}) can be rewritten as
$$y(ch) = y_0 + h\I_\infty(c)^\top\otimes I_m \bfgamma$$
and, consequently, the vector $\bfgamma$ satisfies the equation
\begin{equation}\label{bfgameq}
\bfgamma = \int_0^1 \P_\infty(\tau) \otimes I_m f(y_0+h\I_\infty(\tau)^\top\otimes I_m\bfgamma)\dd \tau,
\end{equation}
with (compare with (\ref{csRKinf})),
\begin{equation}\label{yhnew}
y(h) = y_0+h\gamma_0(y).
\end{equation}

\subsection{Polynomial approximation}\label{ode1_s}

In order to derive a polynomial approximation $\sigma\in\Pi_s$ of (\ref{csRKinf}), it suffices to truncate the infinite series in (\ref{ainf}) after $s$ terms:
\begin{equation}\label{as}
a_s(c,\tau) = \sum_{j=0}^{s-1} \int_0^c P_j(x)\dd x P_j(\tau),
\end{equation}
so that
\begin{equation}\label{a1s}
\dot a_s(c,\tau) = \sum_{j=0}^{s-1} P_j(c) P_j(\tau),
\end{equation}
and, therefore, (\ref{csRKinf}) is approximated by  
\begin{eqnarray}\nonumber
\dot \sigma(ch) &=& \int_0^1\dot a_s(c,\tau)f(\sigma(\tau h))\dd\tau, \qquad \sigma(0)=y_0\in\RR^m,\\ \label{csRKs}
\sigma(ch)        &=& y_0 + h\int_0^1 a_s(c,\tau)f(\sigma(\tau h))\dd\tau, \qquad c\in[0,1],\\ \nonumber
y_1\,:=\,\sigma(h) &=& y_0 + h\int_0^1 f(\sigma(ch))\dd c.
\end{eqnarray} 
\begin{rem}\label{hbvminf} As it was shown in \cite{ABI2019}, the csRK method
\begin{equation}\label{csRK1s}
\begin{array}{c|c} c & a_s(c,\tau)\\ \hline & 1\end{array}, 
\end{equation}
with $a_s(c,\tau)$ given by (\ref{as}), is equivalent to the energy-preserving method, named HBVM$(\infty,s)$, introduced in \cite{BIT2010}.\footnote{In particular, when $s=1$ one obtains the AVF method \cite{QMcL2008}; for a related approach, see also \cite{H2010}.} 
\end{rem}
The following properties hold true \cite{BIT2012_3}:
\begin{equation}\label{y1syh}
y_1-y(h) = O(h^{2s+1}),
\end{equation}
i.e., the approximation procedure has order $2s$. 

\begin{rem}\label{hbvminfcons} If \,$H:\RR^m\rightarrow\RR$,\, and
$$f(y) = J\nabla H(y), \qquad with \qquad J^\top=-J,$$
then \,$H(y_1)=H(y_0)$. In the case of Hamiltonian problems, $H$ is the energy of the system.
Consequently, the csRK method (\ref{csRK1s}) is energy-conserving, as is shown in \cite[Theorem~3]{BIT2012_3}).
\end{rem}

A corresponding vector formulation of the csRK method (\ref{csRK1s}) can be derived by replacing the infinite vectors and matrix in (\ref{PIinf})-(\ref{Xinf}) with
\begin{equation}\label{PIs}
\P_r(c) = \pmatrix{c} P_0(c)\\  \vdots \\ P_{r-1}(c)\endpmatrix, \quad r=s,s+1, \qquad \I_s(c) = \int_0^c \P_s(x)\dd x,
\end{equation}
 and the matrices,
 \begin{equation}\label{Xs}
 \hat{X}_s = \pmatrix{cccc} \xi_0 &-\xi_1\\
                                         \xi_1 & 0 &\ddots\\
                                                  &\ddots &\ddots &-\xi_{s-1}\\
                                                  &            &\xi_{s-1} &0 \\
                                                  \hline
                                                  &            &              &\xi_s
                                                  \endpmatrix \equiv 
                                                  \pmatrix{c} X_s \\
                                                  \hline
                                                        0 \dots 0   \,     \xi_s
                                                  \endpmatrix \in\RR^{s+1\times s}, 
                                              \end{equation}    
such that
\begin{equation}\label{PPPIXs} 
\I_s(c)^\top = P_{s+1}(c)^\top \hat X_s, \quad
\int_0^1 \P_r(\tau)\P_r(\tau)^\top\dd \tau = I_r, \quad \int_0^1 \P_s(\tau)\I_s(\tau)^\top\dd \tau = X_s,
\end{equation}
with $I_r\in\RR^{r\times r}$ the identity matrix. Consequently, with reference to (\ref{as}), one has
\begin{equation}\label{aseq}
a_s(c,\tau) = \I_s(c)^\top\P_s(\tau) = \P_{s+1}(c)^\top \hat X_s \P_s(\tau),
\end{equation}                                              
and, in particular, one may regard the Butcher tableau, equivalent to (\ref{csRK1s})
\begin{equation}\label{csRK1Ws}
\begin{array}{c|c} c & \P_{s+1}(c)^\top \hat X_s \P_s(\tau) \\ \hline & 1\end{array}, 
\end{equation}
as the corresponding $W$-transformation of a HBVM$(\infty,s)$ method \cite{ABI2019}.

Further, according to \cite{BIT2011}, setting the vector
\begin{equation}\label{bfgams}
\bfgamma_s := \pmatrix{c} \gamma_0(\sigma)\\ \vdots \\ \gamma_{s-1}(\sigma)\endpmatrix,
\end{equation} 
with $\gamma_j(\sigma)$ defined as in (\ref{gammaj}), by formally replacing $y$ by $\sigma$, one has that such vector satisfies the equation
(compare with (\ref{bfgameq}))
\begin{equation}\label{bfgameqs}
\bfgamma_s = \int_0^1 \P_s(\tau) \otimes I_m f(y_0+h\I_s(\tau)^\top\otimes I_m\bfgamma_s)\dd \tau,
\end{equation}
with (compare with (\ref{csRKs}) and (\ref{yhnew})),
\begin{equation}\label{y1new}
y_1 = y_0+h\gamma_0(\sigma).
\end{equation}

\section{Approximation of special second-order ODE-IVPs}\label{ode2}

An interesting particular case is that of special second order problems, namely problems in the form 
\begin{equation}\label{ivp2}
\ddot y(t) = f(y(t)), \qquad t\in[0,h], \qquad y(0)=y_0,\,\dot y(0)=\dot y_0\in\RR^m,
\end{equation}
which, in turn, are a special case of the more general problems
\begin{equation}\label{ivp2gen}
\ddot y(t) = f(y(t),\dot y(t)), \qquad t\in[0,h], \qquad y(0)=y_0,\,\dot y(0)=\dot y_0\in\RR^m.
\end{equation}

Let us study, at first, the special problem (\ref{ivp2}), which is very important in many applications (as an example, separable Hamiltonian problems are in such a form), then discussing problem (\ref{ivp2gen}). Setting $\dot y(t) = p(t)$, and expanding the right-hand sides along the Legendre basis gives:
\begin{eqnarray}\label{qp1}
\dot y(ch) &=& \P_\infty(c)^\top\otimes I_m \int_0^1 \P_\infty(\tau)\otimes I_m p(\tau h)\dd\tau,\\ \nonumber
\dot p(ch) &=& \P_\infty(c)^\top\otimes I_m \int_0^1 \P_\infty(\tau)\otimes I_m f(y(\tau h))\dd\tau,\qquad c\in[0,1].
\end{eqnarray}

\begin{rem} It is worth mentioning that, with reference to (\ref{a1c}) and (\ref{ainfeq}), the two previous equations can be rewritten as:
\begin{equation}\label{a2inf}
\dot y(ch) = \int_0^1 \dot a_\infty(c,\tau) p(\tau h)\dd\tau, \qquad
\dot p(ch) = \int_0^1 \dot a_\infty(c,\tau) f(y(\tau h))\dd\tau, \qquad c\in[0,1].
\end{equation}
\end{rem}

Integrating side by side (\ref{qp1}), and imposing the initial condition then gives:
\begin{eqnarray}\nonumber
y(ch) &=& y_0 +h\I_\infty(c)^\top\otimes I_m \int_0^1 \P_\infty(\tau)\otimes I_m p(\tau h)\dd\tau,\\ \label{qp}
p(ch) &=& \dot y_0 + h\I_\infty(c)^\top\otimes I_m \int_0^1 \P_\infty(\tau)\otimes I_m f(y(\tau h))\dd\tau \\
        &=&\dot y_0 + h\I_\infty(c)^\top\otimes I_m \bfgamma ,\qquad c\in[0,1],\nonumber
\end{eqnarray}
with the vector $\bfgamma$ formally still given by (\ref{bgam}).
Substitution of the second equation in (\ref{qp}) into the first one,  taking into account that\,\footnote{In general, $e_i$ will denote the infinite vector whose $j$th entry is $\delta_{ij}$.}
$$\int_0^1\P_\infty(\tau)\dd\tau = e_1 \equiv \pmatrix{c} 1\\ 0\\ \vdots\endpmatrix,\qquad and  
\qquad \I_\infty(c)^\top e_1 = c,$$
and considering (\ref{PPPIXinf}) and again (\ref{bgam}), then gives, 
\begin{eqnarray}\nonumber
y(ch) &=& y_0 +h\I_\infty(c)^\top\otimes I_m \int_0^1 \P_\infty(\tau)\otimes I_m \left[\dot y_0 + h\I_\infty(\tau)^\top\otimes I_m \bfgamma\right]\dd\tau\\
\nonumber
&=& y_0 + ch\dot y_0 + h^2 \I_\infty(c)^\top \otimes I_m \int_0^1 \P_\infty(\tau)\I_\infty(\tau)^\top \dd\tau \otimes I_m \bfgamma \\ \nonumber
&=& y_0 + ch\dot y_0 + h^2 \I_\infty(c)^\top X_\infty \otimes I_m  \bfgamma \\  \label{csRKNy}
&\equiv& y_0 + ch\dot y_0 + h^2\int_0^1 \left[ \I_\infty(c)^\top X_\infty \P_\infty(\tau) \right]\otimes I_m f(y(\tau h))\dd\tau
\end{eqnarray} 
Setting
\begin{equation}\label{csRKN}
\bar a_\infty(c,\tau) =  \I_\infty(c)^\top X_\infty \P_\infty(\tau) \equiv \P_\infty(c)^\top X_\infty^2 \P_\infty(\tau),
\end{equation}
one has that (\ref{qp}) can be rewritten as:
\begin{equation}\label{qch}
y(ch) = y_0 + ch\dot y_0 + h^2\int_0^1 \bar a_\infty(c,\tau) f(y(\tau h))\dd\tau, \qquad c\in[0,1].
\end{equation} 

\begin{rem}\label{pchrem} We observe that, from the second equation in (\ref{qp1}) and (\ref{a2inf}), and considering (\ref{ainfeq}), one derives:
\begin{equation}\label{pch}
p(ch) = \dot y_0 + h \int_0^1 a_\infty(c,\tau) f(y(\tau h))\dd\tau, \qquad c\in[0,1].
\end{equation}
Moreover, it is worth mentioning that, see (\ref{csRKN}), (\ref{ainf}), and (\ref{PPPIXinf}),
\begin{equation}\label{bainfeq}
\bar a_\infty(c,\tau) = \int_0^1 \I_\infty(c)^\top\P_\infty(\xi)\I_\infty(\xi)^\top\P_\infty(\tau)\dd\xi \equiv \int_0^1 a_\infty(c,\xi)a_\infty(\xi,\tau)\dd\xi.
\end{equation}
\end{rem}

From (\ref{qp1}) and (\ref{qch}) one obtains that the values at $h$ will be given by (see (\ref{Xinf})):
\begin{equation}\label{ph}
\dot y(h) ~\equiv~ p(h) ~=~ \dot y_0 + h e_1\top \otimes I_m \bfgamma  \,=\,  \dot y_0 + h \gamma_0(y) \,\equiv\,  \dot y_0 + h\int_0^1 f(y(ch))\dd c,\\
\end{equation} and
\begin{eqnarray}\nonumber
y(h) &=& y_0 + h\dot y_0 + h^2 e_1^\top X_\infty \int_0^1 \P_\infty(c)\otimes I_m f(y(c h)\dd c \\ \nonumber
&=& y_0 + h\dot y_0 + h^2 e_1^\top X_\infty \int_0^1 \P_\infty(c)\otimes I_m f(y(c h)\dd c \\ \nonumber
&=& y_0 + h\dot y_0 + h^2 \int_0^1 (\xi_0 e_1 - \xi_1 e_2)^\top\P_\infty(c)\otimes I_m f(y(c h)\dd c \\ \nonumber
&=& y_0 + h\dot y_0 + h^2  \int_0^1 [\xi_0 P_0(c) - \xi_1 P_1(c)] f(y(c h)\dd c \\ \label{qh}
&\equiv& y_0 + h\dot y_0 + h^2  \int_0^1 (1-c) f(y(c h)\dd c.
\end{eqnarray}
In other words, the exact solution of problem (\ref{ivp2}) is generated by the following csRKN ``method'':
\begin{equation}\label{csRKNtab}
\begin{array}{c|c} c & \bar a_\infty(c,\tau) \\ \hline & 1-c\\ \hline & 1\end{array}\quad \equiv\quad
\begin{array}{c|c} c & \I_\infty(c)^\top X_\infty \P_\infty(\tau) \\ \hline & 1-c\\ \hline & 1\end{array}\quad \equiv\quad
\begin{array}{c|c} c & \P_\infty(c)^\top X_\infty^2 \P_\infty(\tau) \\ \hline & 1-c\\ \hline & 1\end{array}.
\end{equation}
 
\subsection{Vector formulation}\label{ode2_v}

An interesting alternative formulation of the method (\ref{csRKNtab}), akin to (\ref{bfgameq})-(\ref{yhnew}) for first order problems, can be derived by combining (\ref{bgam}) and (\ref{csRKNy}):
\begin{equation}\label{bfgameq2}
\bfgamma = \int_0^1 \P_\infty(\tau)\otimes I_m\, f\hspace{-1mm}\left( y_0 + \tau h\dot y_0 + h^2 \I_\infty(\tau)^\top X_\infty \otimes I_m  \bfgamma\right)\dd\tau,
\end{equation}
with values at $t=h$ given by (compare with (\ref{ph})-(\ref{qh})):
\begin{equation}\label{pqhnew}
\dot y(h) = \dot y_0 + h\gamma_0(y), \qquad y(h) = y_0 + h\dot y_0 + h^2\left( \xi_0 \gamma_0(y)-\xi_1 \gamma_1(y)\right).
\end{equation}

\subsection{The case of the general problem (\ref{ivp2gen})}\label{ode2_gen}

The arguments used above, can be extended to cope with problem (\ref{ivp2gen}) in a straightforward way. In fact, by following similar steps as before, 
(\ref{qch})-(\ref{pch}) now becomes:
\begin{eqnarray}\nonumber
y(ch) &=& y_0 + ch\dot y_0 + h^2  \int_0^1 \bar a_\infty(c,\tau) f(y(\tau h), \dot y(\tau h))\dd\tau,\\
\dot y(ch) &=& \dot y_0 + h  \int_0^1 a_\infty(c,\tau) f(y(\tau h), \dot y(\tau h))\dd\tau, \qquad c\in[0,1], \label{sono2}
\end{eqnarray}
with the values at $h$ given by
\begin{equation}\label{yy1hgen}
\dot y_0 + h  \int_0^1 f(y(c h),\dot y(ch)))\dd c, \qquad y(h) = y_0 +h \dot y_0 + h^2  \int_0^1(1-c) f(y(c h), \dot y(c h))\dd c.
\end{equation}
Consequently, we obtain the general csRKN method
\begin{eqnarray}\nonumber
\lefteqn{
\begin{array}{c|c|c} c & \bar a_\infty(c,\tau) & a_\infty(c,\tau) \\ \hline & 1-c & 1\end{array} ~\equiv~
\begin{array}{c|c|c} c & \I_\infty(c)^\top X_\infty \P_\infty(\tau) & \I_\infty(c)^\top \P_\infty(\tau)  \\ \hline & 1-c & 1\end{array}}
\\ \label{gcsRKNtab}
&\equiv&
\begin{array}{c|c|c} c & \P_\infty(c)^\top X_\infty^2 \P_\infty(\tau) & \P_\infty(c)^\top X_\infty \P_\infty(\tau) \\ \hline & 1-c & 1\end{array},\hspace{3cm}
\end{eqnarray}
in place of (\ref{csRKNtab}).

The bad news is that now the system of equations (\ref{sono2}) has a doubled size, w.r.t. (\ref{qch}). On the other hand, the good news is that, upon modifying the definition of the vector $\bfgamma$ in (\ref{bgam}) as follows,
\begin{equation}\label{bgam1}
\bfgamma := \pmatrix{c} \gamma_0(y,\dot y)\\ \gamma_1(y,\dot y)\\ \vdots\endpmatrix \equiv \int_0^1 \P_\infty(\tau)\otimes I_m f(y(\tau h), \dot y(\tau h))\dd \tau,
\end{equation}
one has that the equations (\ref{sono2}) can be rewritten as
\begin{eqnarray*}\nonumber
y(ch) &=& y_0 + ch\dot y_0 + h^2  \I_\infty(c)^\top X_\infty \otimes I_m \bfgamma,\\
\dot y(ch) &=& \dot y_0 + h \I_\infty(c)^\top \otimes I_m \bfgamma, \qquad c\in[0,1]. %, \label{sonoancora2}
\end{eqnarray*}
Consequently, we obtain again a single equation for the vector $\bfgamma$ defined in (\ref{bgam1}),
\begin{eqnarray}\label{bgam1eq}
\lefteqn{\bfgamma = \int_0^1 \P_\infty(\tau)\otimes I_m }\\ \nonumber  && 
f\hspace{-1mm}\left( y_0 + \tau h\dot y_0 + h^2  \I_\infty(\tau)^\top X_\infty \otimes I_m \bfgamma,\, \dot y_0 + h \I_\infty(c)^\top \otimes I_m \bfgamma\right)\dd \tau.
\end{eqnarray}
Similarly, (\ref{ph})-(\ref{qh}) respectively become:
\begin{equation}\label{y1yhnew}
\dot y(h) = \dot y_0 + h\gamma_0(y,\dot y), \quad y(h) = y_0 + h\dot y_0 + h^2\left( \xi_0 \gamma_0(y,\dot y)-\xi_1 \gamma_1(y,\dot y)\right).
\end{equation}
Remarkably enough, the equations (\ref{bgam1eq})-(\ref{y1yhnew}) are very similar to (\ref{bfgameq2})-(\ref{pqhnew}).

\subsection{Polynomial approximation}\label{ode2s}

Polynomial approximations of degree $s$, $\sigma(ch)\approx y(ch)$ and $\sigma_1(ch)\approx p(ch)$, can be obtained by formal substitution of the matrices $\P_\infty(c), \I_\infty(c)$, and $X_\infty$ in (\ref{qp1})--(\ref{qh}) with the corresponding finite ones, $\P_s(c)$,  $\I_s(c)\equiv \P_{s+1}(c)\hat X_s$, and $X_s$ defined in (\ref{PIs})-(\ref{Xs}). Consequently, following similar steps as above, one obtains:
\begin{eqnarray}\label{uv1}
\dot \sigma(ch) &=& \P_s(c)^\top\otimes I_m \int_0^1 \P_s(\tau)\otimes I_m \sigma_1(\tau h)\dd\tau,\\ \nonumber
\dot \sigma_1(ch) &=& \P_s(c)^\top\otimes I_m \int_0^1 \P_s(\tau)\otimes I_m f(\sigma(\tau h))\dd\tau,\qquad c\in[0,1],
\end{eqnarray}
in place of (\ref{qp1}). Similarly, setting (in place of (\ref{csRKN}))
\begin{equation}\label{csRKNs}
\bar a_s(c,\tau) =  \I_s(c)^\top X_s \P_s(\tau) \equiv \P_{s+1}(c)^\top \hat X_s X_s \P_s(\tau),
\end{equation}
from (\ref{uv1}) one derives:
\begin{equation}\label{uch}
\sigma(ch) = y_0 + ch\dot y_0 + h^2\int_0^1 \bar a_s(c,\tau) f(\sigma(\tau h))\dd\tau, \qquad c\in[0,1],
\end{equation} 
with the new approximations 
$$y_1:= \sigma(h)\qquad and\qquad \dot y_1 := \sigma_1(h)$$ 
given by
\begin{equation}\label{uvh}
\dot y_1 = \dot y_0 + h\int_0^1 f(\sigma(ch))\dd c, \qquad y_1 = y_0 + h\dot y_0 + h^2\int_0^1 (1-c) f(\sigma(ch))\dd c.
\end{equation}

\begin{rem}\label{vchrem} Similarly as (\ref{a2inf}), one derives:
\begin{equation}\label{a2s}
\dot \sigma(ch) = \int_0^1 \dot a_s(c,\tau) \sigma_1(\tau h)\dd\tau, \quad
\dot \sigma_1(ch) = \int_0^1 \dot a_s(c,\tau) f(\sigma(\tau h))\dd\tau, \qquad c\in[0,1].
\end{equation}
Moreover, (compare with (\ref{pch})), from the second equation in (\ref{uv1}), and considering (\ref{aseq}), one obtains:
\begin{equation}\label{vch}
\sigma_1(ch) = \dot y_0 + h \int_0^1 a_s(c,\tau) f(\sigma(\tau h))\dd\tau, \qquad c\in[0,1].
\end{equation}
At last, (compare with (\ref{bainfeq})),  from (\ref{csRKNs}), (\ref{aseq}), and (\ref{PPPIXs}), one has:
\begin{equation}\label{bas}
\bar a_s(c,\tau) = \int_0^1 \I_s(c)^\top\P_s(\xi)\I_s(\xi)^\top\P_s(\tau)\dd\xi \equiv \int_0^1 a_s(c,\xi)a_s(\xi,\tau)\dd\xi.
\end{equation}
\end{rem}

Summing all up, through (\ref{csRKNs})--(\ref{uvh}) we have defined the following csRKN method:
\begin{equation}\label{csRKNstab}
\begin{array}{c|c} c & \bar a_s(c,\tau) \\ \hline & 1-c\\ \hline & 1\end{array}\quad \equiv\quad
\begin{array}{c|c} c & \I_s(c)^\top X_s \P_s(\tau) \\ \hline & 1-c\\ \hline & 1\end{array}\quad \equiv\quad
\begin{array}{c|c} c & \P_{s+1}(c)^\top \hat X_s X_s \P_s(\tau) \\ \hline & 1-c\\ \hline & 1\end{array}.
\end{equation}
\begin{rem}\label{hbvminfs2}
This latter method is equivalent to the csRKN method obtained by applying the HBVM$(\infty,s)$ method to the special second-order problem (\ref{ivp2}) \cite{ABI2019}.
\end{rem}

An alternative formulation of the method (\ref{csRKNstab}) can be obtained by repeating arguments similar to those used  in Section~\ref{ode2_v}. In fact, by setting the vector (compare with (\ref{bfgams})--(\ref{y1new}))
\begin{equation}\label{bfgam2s}
\bfgamma_s := \pmatrix{c} \gamma_j(\sigma)\\ \vdots \\ \gamma_{s-1}(\sigma)\endpmatrix \equiv \int_0^1 \P_s(\tau)\otimes I_m f(\sigma(\tau h))\dd\tau,
\end{equation} 
by virtue of (\ref{csRKNs})-(\ref{uch}) such a vector satisfies the equation
\begin{equation}\label{bfgameq2s}
\bfgamma_s = \int_0^1 \P_s(\tau)\otimes I_m\, f\hspace{-1mm}\left( y_0 + \tau h\dot y_0 + h^2 \I_s(\tau)^\top X_s \otimes I_m  \bfgamma_s\right)\dd\tau,
\end{equation}
with the new approximations given by (compare with (\ref{pqhnew})):
\begin{equation}\label{pqhnew_1}
\dot y_1 = \dot y_0 + h\gamma_0(\sigma), \qquad y_1 = y_0 + h\dot y_0 + h^2\left( \xi_0 \gamma_0(\sigma)-\xi_1 \gamma_1(\sigma)\right).
\end{equation}

For the general problem (\ref{ivp2gen}), following similar steps as above, the generalized csRKN method (\ref{gcsRKNtab}) becomes, with reference to (\ref{aseq}):
\begin{eqnarray}\nonumber
\lefteqn{
\begin{array}{c|c|c} c & \bar a_s(c,\tau) & a_s(c,\tau) \\ \hline & 1-c & 1\end{array} ~\equiv~
\begin{array}{c|c|c} c & \I_s(c)^\top X_s \P_s(\tau) & \I_s(c)^\top \P_s(\tau)  \\ \hline & 1-c & 1\end{array}}
\\ \label{gcsRKNstab-1}
&\equiv&
\begin{array}{c|c|c} c & \P_{s+1}(c)^\top \hat X_s X_s \P_s(\tau) & \P_{s+1}(c)^\top \hat X_s X_s \P_s(\tau) \\ \hline & 1-c & 1\end{array}.\hspace{1cm}
\end{eqnarray}

Also in this case, we can derive a vector formulation, similar to (\ref{bfgam2s})-(\ref{bfgameq2s}). In fact, by defining the vector
\begin{equation}\label{gbfgam2s}
\bfgamma_s := \pmatrix{c} \gamma_j(\sigma,\sigma_1)\\ \vdots \\ \gamma_{s-1}(\sigma,\sigma_1)\endpmatrix \equiv \int_0^1 \P_s(\tau)\otimes I_m f(\sigma(\tau h),\sigma_1(\tau h))\dd\tau,
\end{equation} 
it satisfies the equation
\begin{eqnarray}\label{gbfgameq2s}
\lefteqn{
\bfgamma_s = \int_0^1 \P_s(\tau)\otimes I_m\,}\\ \nonumber
&& f\hspace{-1mm}\left( y_0 + \tau h\dot y_0 + h^2 \I_s(\tau)^\top X_s \otimes I_m  \bfgamma_s,\,
\dot y_0 + h \I_s(\tau)^\top\otimes I_m  \bfgamma_s\right)\dd\tau,
\end{eqnarray}
with the new approximations given by: 
\begin{equation}\label{gpqhnew_1}
\dot y_1 = \dot y_0 + h\gamma_0(\sigma,\sigma_1), \qquad y_1 = y_0 + h\dot y_0 + h^2\left( \xi_0 \gamma_0(\sigma,\sigma_1)-\xi_1 \gamma_1(\sigma,\sigma_1)\right).
\end{equation}

\begin{rem}\label{echecazzo}
By comparing the discrete problems (\ref{bfgameq2s}) and (\ref{gbfgameq2s}), one realizes that they have the same dimension, independently of the fact that the latter one solves the general problem (\ref{ivp2gen}). This fact is even more striking since, as we are going to sketch in the next section, this will be the case for a general $k$th order ODE-IVP.
\end{rem}

\subsection{Approximation of general $k$th-order ODE-IVPs}\label{odek}

Let us consider the case of a general $k$th order problem,
\begin{eqnarray}\nonumber
y^{(k)}(ch) &=& f( y(ch), y^{(1)}(ch),\dots,y^{(k-1)}(ch) ), \qquad c\in[0,1], \\ \label{ivpk}
y^{(i)}(0) &=&y^{(i)}_0\in\RR^m, \qquad i=0,\dots,k-1.
\end{eqnarray}
Repeating similar steps as above, by defining the infinite vector\,\footnote{Hereafter, for sake of brevity, we shall skip the arguments of the Fourier coefficients $\gamma_i$.}
\begin{equation}\label{bgamk}
\bfgamma \equiv \pmatrix{c} \gamma_0\\ \gamma_1\\ \vdots\endpmatrix := \int_0^1 \P_\infty(\tau)\otimes I_m f(y(\tau h), \dots,y^{(k-1)}(\tau h))\dd \tau,
\end{equation}
one has that it satisfies that equation
\begin{eqnarray}\label{bgamkeq}
\lefteqn{\bfgamma = \int_0^1 \P_\infty(\tau)\otimes I_m}\\ \nonumber
&& \hspace{-.5cm}
f\hspace{-1mm}\left( \sum_{i=0}^{k-1} \frac{(\tau h)^i}{i!}y_0^{(i)} +h^k  \I_\infty(\tau)^\top X_\infty^{k-1} \otimes I_m \bfgamma,\dots,y_0^{(k-1)} + h \I_\infty(\tau)^\top \otimes I_m \bfgamma\right)\dd \tau,
\end{eqnarray}
with the values at $t=h$ given by:
\begin{equation}\label{y1k}
y^{(i)}(h) = \sum_{j=0}^{k-1-i} \frac{h^j}{j!} y^{(j+i)}_0 + h^{k-i}\sum_{j=0}^{k-1-i} b_j^{(k-1-i)}\gamma_j, \qquad i=0,\dots,k-1,
\end{equation}
$b_j^{(k-1-i)}$,\, $j=0,\dots,k-1-i$,\, being the  $(j+1)$st entry on the first row of the matrix $$X_\infty^{k-1-i}, \qquad i=0,\dots,k-1.$$

A polynomial approximation of degree $s$ (resulting, as usual, into an order $2s$ method) can be derived by formally substituting, in the equation  (\ref{bgamkeq}), $\P_\infty(\tau)$ and $\I_\infty(\tau)$ with $\P_s(\tau)$ and $\I_s(\tau)$, respectively,\footnote{In so doing, the vector $\bfgamma$ now belongs to $\RR^{sm}$.} with the new approximations $y_1^{(i)}\approx y^{(i)}(h)$, given by:
\begin{equation}\label{y1ks}
y^{(i)}_1 \,:=\, \sum_{j=0}^{k-1-i} \frac{h^j}{j!} y^{(j+i)}_0 + h^{k-i}\sum_{j=0}^{k-1-i} b_j^{(k-1-i)}\gamma_j, \qquad i=0,\dots,k-1,
\end{equation}
$b_j^{(k-1-i)}$,\, $j=0,\dots,k-1-i$,\, being now the $(j+1)$st entry on the first row of the matrix $$X_s^{k-1-i},  \qquad i=0,\dots,k-1.$$

\section{Discretization}\label{discre}

\begin{flushright}
%``As is well known, even many relatively simple integrals\\ cannot be expressed in finite terms of elementary functions,\\ and thus must be evaluated by numerical methods.'' \\ Dahlquist and Björk \cite[page.~521]{DaBj2008}. 
{\em
``As is well known, even many relatively simple\\ integrals cannot be expressed in finite terms of\\ elementary functions, and thus must be evalu-\\ ated by numerical methods.'' }\mbox{\hspace{2.55cm}} \\ Dahlquist and Björk \cite[page.~521]{DaBj2008}
\end{flushright}
As is clear, the integrals involved in the Fourier coefficients (\ref{bfgameqs}) need to be approximated by using a suitable quadrature rule, which we choose as the interpolatory Gauss-Legendre quadrature of order $2k$, with $k\ge s$, having weights and abscissae $(b_i,c_i)$, $i=1,\dots,k$. In so doing the vector of the Fourier coefficients (\ref{bfgams})-(\ref{bfgameqs}) becomes
\begin{equation}\label{hbfgams}
\hat{\bfgamma}_s := \pmatrix{c} \hat\gamma_0\\ \vdots \\ \hat\gamma_{s-1}\endpmatrix,
\end{equation} 
satisfying the equation (compare with (\ref{bfgameqs}))
\begin{equation}\label{hbfgameqs}
\hat{\bfgamma}_s = \sum_{j=1}^k b_j \P_s(c_j) \otimes I_m f(y_0+h\I_s(c_j))\otimes I_m\hat{\bfgamma}_s)
\end{equation}
with (compare with (\ref{y1new})),
\begin{equation}\label{hy1new}
y_1 = y_0+h\hat\gamma_0.
\end{equation}
We observe that, by introducing the matrices (see (\ref{PIs}))
\begin{eqnarray}\nonumber
&&\P_r := \pmatrix{c} \P_r(c_1)^\top \\ \vdots \\ \P_r(c_k)^\top \endpmatrix\in\RR^{k\times r},\quad
\I_s := \pmatrix{c} \I_s(c_1)^\top \\ \vdots \\ \I_s(c_k)^\top \endpmatrix\in\RR^{k\times s},\\ \label{PIO}
&&\Omega = \pmatrix{ccc} b_1\\ &\ddots\\ &&b_k\endpmatrix,
\end{eqnarray} 
one has, with reference to (\ref{Xs}),
\begin{equation}\label{PIOprop}
\I_s = \P_{s+1}\hat X_s, \qquad \P_s^\top\Omega \P_s = I_s, \qquad  \P_s^\top\Omega \I_s = X_s.
\end{equation}
Further, by also introducing the vector ~$\bfe = (1,\dots,1)^\top\in\RR^k$,~ one obtains that (\ref{hbfgameqs}) can be rewritten as
\begin{equation}\label{hbfgams1}
\hat{\bfgamma}_s = \P_s^\top\Omega\otimes I_m f(\bfe\otimes y_0 + h\I_s\otimes I_m \hat{\bfgamma}_s).
\end{equation}
Setting
\begin{equation}\label{Yi}
Y \equiv \pmatrix{c} Y_1\\ \vdots \\ Y_k\endpmatrix := \bfe\otimes y_0 + h\I_s\otimes I_m \hat{\bfgamma}_s
\end{equation}
the vector of the stages of the corresponding Runge-Kutta method, from (\ref{hbfgams1})-(\ref{Yi}) one obtains the stage equation
\begin{equation}\label{Yeq}
Y = \bfe\otimes y_0 + h\I_s\P_s^\top\Omega\otimes I_m f(Y),
\end{equation}
with the new approximation
\begin{equation}\label{y1eq}
y_1 = y_0 + h\sum_{i=1}^k b_i f(Y_i).
\end{equation}
Summing all up, (\ref{hbfgams1})--(\ref{y1eq}) define the $k$-stage Runge-Kutta method
\begin{equation}\label{HBVMtab}
\begin{array}{c|c} \bfc & \I_s\P_s^\top\Omega \\ \hline & \bfb^\top\end{array} ~,\qquad \bfb = (b_1,\dots,b_k)^\top, \qquad \bfc = (c_1,\dots,c_k)^\top,
\end{equation}
named HBVM$(k,s)$ \cite{LIMbook2016,BI2018ax,BIT2010}. It is worth mentioning that the Butcher matrix of the method, with reference to the coefficients of the csRK (\ref{aseq}), is given by:
$$  \I_s\P_s^\top\Omega \equiv \pmatrix{c} b_ja_s(c_i,c_j)\endpmatrix \in\RR^{k\times k}.$$
In particular, when $k=s$ one obtains the $s$-stage Gauss collocation method. However, the use of values $k>s$ (and even $k\gg s$) is useful, in view of deriving energy-conserving methods for Hamiltonian systems \cite{BIT2010,BIT2012_3,LIMbook2016,BI2018ax}.

\begin{rem}\label{emegliogamma}
As is clear, the formulation (\ref{hbfgameqs})-(\ref{hy1new}) is computationally more effective than (\ref{Yeq})-(\ref{y1eq}), having the former (block) dimension $s$, independently of $k$, which is the dimension of the latter formulation. As observed in \cite{BIT2011}, this allows the use of relatively large values of $k$, without increasing too much the computational cost.
\end{rem}

Similar arguments can be repeated in the case of the polynomial approximations for problems (\ref{ivp2}), (\ref{ivp2gen}), and (\ref{ivpk}): we here sketch only those concerning the csRKN (\ref{csRKNstab}), providing the correct implementation for the HBVM$(k,s)$ method for the special second-order problem (\ref{ivp2}). By formally using the same approximate Fourier coefficients (\ref{hbfgams})  in place of (\ref{bfgam2s}), one has that (\ref{bfgameq2s}) is replaced by the following discrete counterpart,
\begin{equation}\label{hbfgameq2s}
\hat{\bfgamma}_s = \sum_{j=1}^k b_j \P_s(c_j) \otimes I_m f(y_0 +c_jh \dot y_0+h^2\I_s(c_j))X_s\otimes I_m\hat{\bfgamma}_s)
\end{equation}
with the new approximations (compare with (\ref{pqhnew_1})) given by
\begin{equation}\label{hpqhnew_1}
\dot y_1 = \dot y_0 + h\hat\gamma_0, \qquad y_1 = y_0 + h\dot y_0 + h^2\left( \xi_0 \hat\gamma_0-\xi_1 \hat\gamma_1\right).
\end{equation}
Similarly as in the first order case, (\ref{hbfgameq2s}) can be rewritten as
\begin{equation}\label{hbfgams2}
\hat{\bfgamma}_s = \P_s^\top\Omega\otimes I_m f(\bfe\otimes y_0 + h\bfc\otimes \dot y_0+ h^2\I_sX_s\otimes I_m \hat{\bfgamma}_s),
\end{equation}
where $\bfc$ is the vector of the abscissae. Again, the vector
\begin{equation}\label{Yi2}
Y \equiv \pmatrix{c} Y_1\\ \vdots \\ Y_k\endpmatrix := \bfe\otimes y_0 + h\bfc\otimes \dot y_0+ h^2\I_sX_s\otimes I_m \hat{\bfgamma}_s
\end{equation}
is the stage vector of a $k$-stage RKN method. In particular, from (\ref{hbfgams2}) and (\ref{Yi2}) one obtains
\begin{equation}\label{chitelofafare}
Y =  \bfe\otimes y_0 + h\bfc\otimes \dot y_0+ h^2\I_sX_s \P_s^\top\Omega\otimes I_m f(Y),
\end{equation}
and some algebra shows that the new approximations are given by
\begin{equation}\label{chitelofafare_1}
\dot y_1 = \dot y_0 + h\sum_{i=1}^k b_if(Y_i), \qquad y_1 = y_0 + h\dot y_0 + h^2\sum_{i=1}^k b_i(1-c_i)f(Y_i).
\end{equation}
Consequently, we are speaking about the following RKN method:\footnote{Here, $\circ$ denotes the Hadamard, i.e. componentwise, product.}
\begin{equation}\label{HBVM2tab}
\begin{array}{c|c} \bfc & \I_sX_s\P_s^\top\Omega \\ \hline & [\bfb\circ(\bfe-\bfc)]^\top \\ \hline & \bfb^\top\end{array}.
\end{equation}
The Butcher tableau (\ref{HBVM2tab}) defines the RKN formulation of a HBVM$(k,s)$ method \cite{ABI2019,BIT2011,LIMbook2016}. Is is worth mentioning that, with reference to (\ref{csRKNs}) and (\ref{bas})
$$\I_sX_s\P_s^\top\Omega  \equiv \pmatrix{c} b_j \bar a_s(c_i,c_j)\endpmatrix\equiv \pmatrix{c} b_j\int_0^1 a_s(c_i,\tau)a_s(\tau,c_j)\dd\tau\endpmatrix \in\RR^{k\times k}.$$

\begin{rem}As observed in Remark~\ref{emegliogamma}, also in this case, in consideration that $k\ge s$, the formulation (\ref{hpqhnew_1})-(\ref{hbfgams2}) of the method is much more efficient than the usual one, given by (\ref{chitelofafare})-(\ref{chitelofafare_1}), also in view of the use of values of $k\gg s$, as in the case of separable Hamiltonian problems. \end{rem}

We conclude this section by recalling that:
\begin{itemize}
\item very efficient iterative procedures exist for solving the discrete problems (\ref{hbfgams1}) and (\ref{hbfgams2}) generated by a HBVM method (see, e.g., \cite{BIT2011});\footnote{We mention that a state-of-art Matlab code is available at the website of the book \cite{LIMbook2016}.}

\bigskip
\item when HBVMs are used as {\em spectral methods in time}, i.e., choosing values of $s$ and $k>s$ so that no accuracy improvement can be obtained, for the considered finite-precision arithmetic and  timestep $h$ used \cite{ABI2020,BIMR2019,BMR2019}, then there is no practical difference between the discrete methods and their continuous-stage counterparts.
 
\end{itemize}

\section{Conclusions}\label{fine}

In this paper we have recalled the basic facts, also reporting some new insight, on the energy-conserving class of Runge-Kutta methods named HBVMs. Such method have been here studied within the framework of continuous-stage Runge-Kutta methods. The extension to second-order  problems has been also studied, providing a natural continuous-stage Runge-Kutta-Nystr\"om formulation of the methods. Further, also the extension to general $k$th-order problems has been sketched. The relation with the fully discrete methods has been also recalled, thus showing the usefulness of using such a framework to study the fully discrete methods.

\vspace{6pt}

\paragraph{Conflict of interests.} The authors declare no conflict of interest.

%%%%%%%%%%%%%%%%%%%%%%%%%%%%%%%%%%%%%%%%%%
\end{document}